\documentstyle[10pt, amscd]{amsart}
\renewcommand{\epsilon}{\varepsilon}
\newcommand{\newsection}[1]
{\subsection{#1}\setcounter{theorem}{0}\setcounter{equation}{0}
\par\noindent}

\newtheorem{theorem}{Theorem}

\newtheorem{lemma}[theorem]{Lemma}
\newtheorem{corr}[theorem]{Corollary}
\newtheorem{prop}[theorem]{Prop}

\newtheorem{deff}[theorem]{Definition}

\newcommand{\bth}{\begin{theorem}}
\newcommand{\ble}{\begin{lemma}}
\newcommand{\bco}{\begin{corr}}
\newcommand{\bdeff}{\begin{deff}}
\newcommand{\bprop}{\begin{prop}}
\newcommand{\eth}{\end{theorem}}
\newcommand{\ele}{\end{lemma}}
\newcommand{\eco}{\end{corr}}
\newcommand{\edeff}{\end{deff}}
\newcommand{\eprop}{\end{prop}}

\newcommand{\Rn}{{\Bbb R}^n}

\newcommand{\supp}{\text{supp }}
\renewcommand{\Pi}{\varPi}

\renewcommand{\epsilon}{\varepsilon}

\newcommand{\Rt}{{\Bbb R}^2}
\newcommand{\Mdel}{{{\cal M}}^\delta}
\newcommand{\dist}{{\text{dist}}}

\begin{document}

\title[Negative results for Nikodym maximal functions]
{Negative results for Nikodym maximal functions and related oscillatory
integrals in curved space}
\thanks{The first author was supported in part by an NSF postdoctoral fellowship.
 The second author was supported in part by the NSF and was on leave from
 UCLA}
\author{William P. Minicozzi II}
\author{Christopher D. Sogge}
\address{Department of Mathematics, The Johns Hopkins University, Baltimore,
Maryland 21218}

\maketitle

\newsection{Introduction}

In 1972 Carleson and Sj\"olin \cite{CS} proved an optimal theorem for
spherical summation operators in the plane.  Specifically, they showed
that the Fourier multiplier operators corresponding to
$m_\delta(\xi)=(1-|\xi|)^\delta_+$ are bounded on $L^p({\Bbb R}^2)$,
$p\ge 4$ if $\delta>\delta(p)=2(1/2-1/p)-1/2$.  Since the kernel of
this summation operator (the inverse Fourier transform of $m_\delta$)
behaves at infinity like $\sum e^{\pm i|x|}/|x|^{3/2+\delta}$, they
obtained this result by proving the essentially equivalent theorem that
\begin{equation}\label{slam}
S_\lambda f(x)=\int e^{i\lambda|x-y|} a(x,y)f(y)\, dy
\end{equation}
satisfies 
\begin{equation}\label{csbound}
\|S_\lambda f\|_{L^4({\Bbb R}^2)}\le C_\epsilon
\lambda^{-1/2+\epsilon}\|f\|_{L^4({\Bbb R}^2)}, \, \, \, \lambda\ge 1,
\varepsilon>0,
\end{equation}
if $a\in C^\infty_0(\Rt \times \Rt)$ vanishes near the diagonal where
$x=y$.  Using a scaling argument, one finds that this yields the preceding
multiplier theorem when $p=4$, and the other cases follow from
interpolating with the easy estimate corresponding to $p=\infty$.

Carleson and Sj\"olin actually proved a stronger result.  They
considered oscillatory integral operators of the form
\begin{equation}\label{cs}
T_\lambda f(x)=\int e^{i\lambda \phi(x,t)}a(x,t) f(t)\, dt,
\end{equation}
where now $a, \phi\in C^\infty(\Rt \times{\Bbb R})$ and moreover the
real phase function is assumed to satisfy the Carleson-Sj\"olin
condition that
\begin{equation}\label{cscond}
\text{det }
\begin{pmatrix}
\phi''_{x_1t} & \phi'''_{x_1tt}
\\
\phi''_{x_2t} & \phi'''_{x_2tt}
\end{pmatrix}\ne 0,
\, \, \, \text{on  supp }a.
\end{equation}
Under these hypotheses they proved the following stronger more general
version of \eqref{csbound}:
\begin{equation}\label{csbound2}
\|T_\lambda f\|_{L^4(\Rt)}\le C_\epsilon
\lambda^{-1/2+\epsilon}\|f\|_{L^4({\Bbb R})}, \, \, \epsilon>0.
\end{equation}

In the other direction Fefferman \cite{F1} had earlier showed that the
multiplier operators corresponding to $\delta=0$, that is, the ball
multiplier operators with $m_0(\xi)=\chi_{|\xi|\le1}$ are never bounded
on $L^p(\Rn)$ if $n\ge2$ and $p\ne2$.  The proof in this seminal paper
involved using Besicovitch's construction that there are sets in the plane of
measure zero containing a unit line segment in every direction.  Using
related ideas, in \cite{F2}, Fefferman was able to give an independent
proof of the Carleson-Sj\"olin multiplier theorem which had a more
geometric flavor.  Many of the recent results in the subject
use ideas from Fefferman's work.

Following \cite{F2} in part, C{\'o}rdoba \cite{Co} gave another proof of the
Carleson-Sj\"olin theorem.  Using a straightforward orthogonality
argument which exploited the fact that the critical estimate involves
$L^4$ and $4=2\cdot2$, C{\'o}rdoba showed that the multiplier theorem
follows from optimal bounds for the ``Nikodym maximal operators" in the
plane.  Specifically, if $T^\delta$ denotes a $\delta$-neighborhood of
a unit line segment in $\Rt$ and if 
\begin{equation}\label{nik}
({{{\cal M}}}^\delta f)(x)=\sup_{x\in T^\delta}
|T^\delta|^{-1}\int_{T^\delta}
|f(y)|\, dy,
\end{equation}
C{\'o}rdoba showed that when $\varepsilon>0$ and $0<\delta\le 1$,
\begin{equation}\label{Cordoba}
\|{\cal M}^\delta f\|_{L^2(\Rt)}\le C_\epsilon \delta^{-\epsilon}
\|f\|_{L^2(\Rt)}.
\end{equation}
C{\'o}rdoba also conjectured that for higher dimensions one should have the
optimal bounds
\begin{equation}\label{cordconj}
\|{\cal M}^\delta f\|_{L^q(\Rn)}\le C_{p,\epsilon}
\delta^{1-n/p-\epsilon}\|f\|_{L^p(\Rn)}, \,
\,
q=(n-1)p', \, \, 1\le p\le n,
\end{equation}
assuming as before that $0<\delta\le 1$ and $\epsilon>0$.  Here, and in
what follows, $p'=p/(p-1)$ denotes the exponent which is conjugate to $p$.

While this estimate is not known there are many partial results.  First
of all Christ, Duandikoetxea and Rubio de Francia \cite{CDF} showed
that \eqref{cordconj} holds when $p\le (n+1)/2$.  (See also Drury
\cite{D} for related estimates.)  This estimate then was improved in an
important paper of Bourgain \cite{B1}, in which it was shown that when
$n\ge3$ \eqref{cordconj} a slightly weaker version of \eqref{cordconj}
(with other norms in the left)
holds for certain $(n+1)/2<p\le p_n$, where $p_n$ is
given by a certain recursive relation arising from an induction
argument on the dimension $n$.  Wolff \cite{W} then improved Bourgain's
result, showing that when $n\ge3$ \eqref{cordconj} holds for $p\le
(n+2)/2$.

In this paper we shall show how an argument of Bourgain \cite{B1} and
Wolff \cite{W} can be used to show that on a Riemannian manifold of dimension
$n$ an analog of \eqref{cordconj} holds for $p\le (n+1)/2$, if in 
\eqref{nik} $T^\delta$ are $\delta$-neighborhoods of geodesics of an
appropriate length and the norms are defined using
the volume element.  In odd dimensions we shall show that this result
is optimal.  Specifically, we shall provide an example of a Riemannian
manifold for which the analog of \eqref{cordconj} does not hold for
any $p>[(n+2)/2]$, if $[(n+2)/2]$ denotes the greatest integer $\le
(n+2)/2$.  We do this by showing that in curved space Nikodym-type sets
of dimension $[(n+2)/2]$ may exist.  The
aforementioned positive results for ${\cal M}^\delta$ imply that such sets
must always have dimension $\ge (n+1)/2$.  The Nikodym-type sets we
construct turn out to be smooth submanifolds and since $(n+1)/2$ is a
half integer for even $n$, this explains the gap between the negative
and positive results for the general case here.  Similar
numerology also
arose in some negative results of Bourgain \cite{B2} for oscillatory
integrals.

The main idea behind our constructions comes from the proof of positive
results for the Euclidean setting of Bourgain \cite{B1} and Wolff
\cite{W}.  In each of these papers a key step involves reducing to
estimates for ${\cal M}^\delta$ involving lower dimensions $2\le m<n$.
To extend these proofs in a trivial way to a curved space setting one would
need that there are many totally geodesic submanifolds of dimension
$m$.  Unfortunately, for non-Euclidean manifolds, it is of course rare to have this
if $m\ne 1$ or $n$, and all of our counterexamples are built around this
fact.  On the other hand, we should point out that our results suggest
that the worst cases for \eqref{cordconj} and the related oscillatory
integral estimates described below might involve metrics whose
sectional curvatures degenerate to high order along lower dimensional sets.

Let us now turn to the related negative results for oscillatory
integrals.  To put them in context, we first need to recall a work of
H\"ormander \cite{H1}.  In this paper, the proof of Carleson-Sj\"olin
\cite{CS} was simplified and H\"ormander improved their oscillatory
integral estimate \eqref{csbound2} by showing that
\begin{equation}\label{csbound3}
\|T_\lambda f\|_{L^q(\Rt)}\le C_q\lambda^{-2/q}\|f\|_{L^p({\Bbb R})}, \,
\, \, 4<q\le \infty, \, \, p=3p'.
\end{equation}
This result can be
seen to be best possible.  H\"ormander also formulated a natural
extension of the Carleson-Sj\"olin condition for real phase functions
$\phi(x,t)\in C^\infty(\Rn\times {\Bbb R}^{n-1})$ and raised the
problem of trying to generalize \eqref{csbound3} to higher dimensions.
This higher dimensional version of the Carleson-Sj\"olin condition
\eqref{cscond} can be formulated as follows.  First one requires that the
mixed Hessian of the phase function have maximal rank on $\text{supp
}a$, that is,
\begin{equation}\label{hor1}
\text{rank }(\partial^2\phi/\partial x_j\partial t_k)\equiv n-1.
\end{equation}
If this condition is met and if we fix $x=x_0\in \text{supp}_x\, a$, then
\begin{equation}\label{surface}
\Sigma_{x_0}=\{\nabla_x\phi(x_0,t): \, \, t\in {\cal N}\}
\end{equation}
is a smooth (immersed) hypersurface in $\Rn$ if ${\cal N}$ is a small
neighborhood of $\{t: \, a(x_0,t)\ne 0\}$.  The other part of the
Carleson-Sj\"olin condition is that
\begin{equation}\label{hor2}
h_{jk} \, \text{is nondegenerate on } \, \Sigma_{x_0},
\end{equation}
if $h_{jk}$ denotes the second fundamental form of $\Sigma_{x_0}$
induced by the Euclidean metric on $\Rn$.  These conditions are easily seen
to be invariant and it is clear that they are equivalent to
\eqref{cscond} when $n=2$.  Assuming them, H\"ormander asked whether
bounds of the form
\begin{equation}\label{horconj}
\|T_\lambda f\|_{L^q(\Rn)}\le C_q\lambda^{-n/q}
\|f\|_{L^p({\Bbb R}^{n-1})}, \, \, 2n/(n-1)<q\le \infty, \, \, q=(n+1)p'/(n-1)
\end{equation}
hold when $n\ge3$.

The first general result of this type is due to Stein
\cite{St2} who showed that when $n\ge3$, \eqref{horconj} holds for
$q\ge 2(n+1)/(n-1)$, generalizing the earlier $L^2$ restriction theorem
of Stein and Tomas \cite{T}.  In the other direction, Bourgain
\cite{B1} provided a striking example showing how, at least for odd
$n$, Stein's result is optimal.  When $n=3$, following Stein
\cite{St3}, it is particularly easy to describe Bourgain's example.
One simply takes
\begin{equation}\label{bad}
\phi(x,t)=x_1t_1+x_2t_2+<A(x_3)t,t>,
\end{equation}
where, say,
$$A(x_3)=
\begin{pmatrix}
1&x_3
\\
x_3&x_3^2
\end{pmatrix},
$$
so that
$$\text{rank }A\equiv 1, \, \, \text{but } \, \text{rank }A'=2.$$
Clearly, \eqref{hor1} holds and since $A'$ has full rank the other
part, \eqref{hor2}, of the Carleson-Sj\"olin condition must hold.
Since $\text{rank }\phi''_{tt}\equiv1$ one can use stationary phase to
see that if the amplitude $a$ of $T_\lambda$ is nonnegative and if a fixed
$f\in C^\infty_0$ equals one
on $\text{supp}_t a\ne\emptyset$, then $|T_\lambda f(x)|\approx
\lambda^{-1/2}$
for large $\lambda>1$, if $x$ is a distance $O(\lambda^{-1})$ from
$\text{supp}_x a \cap \{(x',x_3): \, x'\in \text{range }A(x_3)\}$.
Hence, $\|T_\lambda f\|_q/\|f\|_\infty \ge C\lambda^{-1/2-1/q}$,
showing that \eqref{horconj} cannot hold here when $q<4$, as claimed.

The mechanism behind this example that $\text{rank }\phi''_{tt}<n-1$
everywhere does not seem possible if, unlike the preceding case, the
second fundamental forms in the second part of the Carleson-Sj\"olin
condition are always positive definite.  The latter happens in the
model case where $\phi(x,t)$ is the Riemannian distance between $x$ and
$t$ with $t$ belonging to an appropriate hypersurface and $x$ belonging
to the compliment.  In this case, the second fundamental forms cannot
have positive signature since, by Gauss' lemma, the surfaces
\eqref{surface} are just the cospheres $\{\xi: \, \sum_{j,k=1}^n
g^{jk}(x_0)\xi_j\xi_k =1\}$, with $g^{jk}=(g_{jk})^{-1}$ denoting the
cometric coming from the Riemannian metric $\sum g_{jk}dx_jdx_k$ on the
manifold $M^n$.

Because of this one might hope for better results for $T_\lambda$ if,
as above, one considers the model case where the phase functions come from
a Riemannian metric.  Here too, though, things may break down.  Indeed,
using the same counterexamples for \eqref{cordconj}, we shall show
that, even if one considers weaker estimates involving now
\begin{equation}
\label{osc}
S_\lambda f(x)=\int_{M^n}e^{i\lambda \text{dist}(x,y)}a(x,y)f(y)\, dy,
\end{equation}
then 
\begin{multline}\label{oscest}
\|S_\lambda f\|_{L^p(M^n)} \le
C_{q,\varepsilon}\lambda^{-n/q+\varepsilon} \|f\|_{L^q(M^n)},
\\
2n/(n-1)<q\le \infty, \, \, q=(n+1)p'/(n-1), \, \, \varepsilon>0,
\end{multline}
need not hold for $n=3$ if $3<q<10/3$.  Here,
$\text{dist}(\cdot,\cdot)$ is the distance coming from the metric
$g_{jk}$ on $M^n$, and, as before, the amplitude is assumed to be
$C^\infty_0$ and to vanish near the diagonal to insure that the phase
function is smooth.  In this context, we sharpen a negative result of
Bourgain \cite{B2} who showed that \eqref{horconj} generically breaks
down if $q<118/39$.  As with the Nikodym maximal functions the metrics
can be taken to be real analytic and arbitrarily close to the Euclidean
one.  The constructions also give negative results for $n>3$.

\newsection{Negative results for the Nikodym maximal function when
$n=3$}

Before  focusing on the three-dimensional case, let us describe the
general setup.  Let $M^n$ be a complete $n$-dimensional Riemannian
manifold.  We shall consider all geodesics $\gamma_x$ containing a
given point $x\in M^n$ of length $|\gamma_x|=r$.  We then for
$0<\delta\le 1$ let $T^\delta_{\gamma_x}$ denote a tubular neighborhood
of width $\delta$ around $\gamma_x$ and define
\begin{equation}\label{2.1}
\Mdel f(x)=\sup_{x\in\gamma_x, \,
|\gamma_x|=r}|T^\delta_{\gamma_x}|^{-1}
\int_{T^\delta_{\gamma_x}}|f(y)|\, dy.
\end{equation}
If we then fix a compact subset $K\subset M^n$, we shall be concerned
with the problem of deciding when bounds of the form
\begin{equation}\label{2.2}
\|\Mdel f\|_{L^q(K)}\le
C_{p,\varepsilon}\delta^{1-n/p-\varepsilon}\|f\|_{L^p}, \, \,
q=(n-1)p', \, \, \varepsilon>0, \, \, \text{supp }f\subset K
\end{equation}
can hold, assuming of course that $1\le p\le n$.  Later we shall give a
simple argument based on \cite{B1} and \cite{W} showing that if $r$ as
above is small enough then the analog of the Euclidean results in
\cite{CDF} always hold.  Specifically, we shall see that \eqref{2.2}
holds on an arbitrary manifold if $1\le p\le (n+1)/2$.  Before doing
this, we shall show that for odd dimensions this result is sharp in the
sense that there are odd-dimensional manifolds for which \eqref{2.2}
cannot hold for any $p>(n+1)/2$ regardless of how small we choose the
fixed number $r$ to be.  For even $n$ we shall show that \eqref{2.2}
breaks down for $p>(n+2)/2$.  We shall also give a simple explanation
of the difference between even and odd dimensions for our type of
constructions. 

Let us start out with the negative results for Nikodym maximal functions
when $n=3$ since this is the simplest case.  Here we wish to show that
\eqref{2.2} need not hold on a given curved three-dimensional
Riemannian manifold if $p>2$.  The main step involves the following
simple lemma.

\ble\label{mainlemma}
Let $\alpha\in C^\infty(\Bbb R)$ satisfy $-1<\alpha<1$ and
$\alpha(0)=0$ and set $\alpha^{(-1)}(t)=\int_0^t\alpha(s)\, ds$.  Let
\begin{equation}\label{2.3}
p(x,\xi)=\sqrt{|\xi|^2+2\alpha(x_2)\xi_1\xi_3}
\end{equation}
be the symbol of the cometric
$\sum g^{jk}(x)d\xi_jd\xi_k=d\xi^2+2\alpha(x_2)d\xi_1d\xi_3$ on
$T^*\Bbb R^3$.  Then for fixed $x_1\in \Bbb R$, and
$-\pi/2<\theta<\pi/2$
\begin{equation}\label{2.4}
t\to x(x_1,\theta;t)=(\, x_1+t\sin\theta, \, t\cos\theta, \, \sin\theta
\alpha^{(-1)}(t\cos\theta)/\cos\theta\, )
\end{equation}
is a geodesic for the corresponding metric $\sum g_{jk}(x)dx_jdx_k$ on
$T\Bbb R^3$, where $g_{jk}=(g^{jk})^{-1}$.  Furthermore, the Jacobian
of the map
\begin{equation}\label{2.5}
(x_1,\theta,t)\to x(x_1,\theta;t)
\end{equation}
equals $|\alpha^{(-1)}(t)|$ when $\theta=0$.
\ele

\begin{pf}
The last assertion involves a straightforward calculation.  To verify
that the curves \eqref{2.4} are geodesics for our metric, we need to
recall that if $(x(t),\xi(t))$ satisfies Hamilton's equation
\begin{equation}\label{2.6}
dx/dt=\partial p/\partial\xi, \, \, d\xi/dt=-\partial p/\partial x,
\end{equation}
then $t\to x(t)$ is geodesic.  (See, e.g., Appendix C in \cite{H2}.)
Furthermore, since $p$ must be constant on its integral curves, if we
take
$$x(0)=(x_1,0,0), \, \, \xi(0)=(\sin\theta,\cos\theta,0)$$
as initial conditions, then, since $p(x(0),\xi(0))=1$, \eqref{2.6}
becomes in our case
$$dx/dt=(\xi_1+\alpha(x_2)\xi_3,\xi_2,\xi_3+\alpha(x_2)\xi_1), \, \,
d\xi/dt=-(0,\alpha'(x_2)\xi_1\xi_3,0).$$
Our initial condition then yields
$\xi(t)=\xi(0)=(\sin\theta,\cos\theta,0)$.  If we plug this into the
formula for $dx/dt$ we conclude that
$(x_1(t),x_2(t))=(x_1+t\sin\theta,t\cos\theta)$, as desired. We then
integrate the last variable to obtain
$$x_3(t)=\int_0^t\sin\theta \, \alpha(s\cos\theta)\, ds,$$
yielding the remaining part of \eqref{2.4}
\end{pf}

To apply the lemma take
\begin{equation}\label{alpha}
\alpha(s)=e^{1/s}, s<0, \, \, \text{and } \, \alpha(s)=0, \, \, s\ge 0,
\end{equation}
and let $\sum g_{jk}dx_jdx_k$ be the metric corresponding to the cometric
$d\xi^2+2\alpha(x_2)d\xi_1d\xi_3$.  The metric then agrees with the
Euclidean one for $x_2\ge0$.  Moreover, since $\alpha^{(-1)}(s)=0$ for
$s\ge0$, the lemma implies that there is an open neighborhood ${\cal
N}\subset \{x\in {\Bbb R}^3: \, x_2<0\}$ of the half-axis where $x_2<0$,
$x_1=x_3=0$ so that if $x\in {\cal N}$ there is a unique geodesic
$\gamma_x$ containing $x$ and having the property that when $x_2\ge0$
$\gamma_x$ is contained in the two-plane $x_3=0$.  If we then, for a
given $c>0$, let
$$f_\delta(x)=1 \, \text{if } \, x_2>0, \, |(x_1,x_2)|<c\, \, \text{and
} \, |x_3|<\delta, \, \text{and } f_\delta(x)=0 \, \, \text{otherwise},$$
it follows that for small fixed $x_2<0$, $\Mdel f_\delta(x)$ must be bounded
from below by a positive constant on some nonempty Euclidean ball $B$
centered at $(0,x_2,0)$.  Hence,
$$\|\Mdel f_\delta\|_{L^1(B)}\, /\, \|f_\delta\|_{L^p} \, \ge
c_0\delta^{-1/p}$$
for some $c_0>0$ depending on $B$ and $c>0$ above.  Since
$$3/p-1<1/p \, \, \, \text{when } \, p>2,$$
we conclude that \eqref{2.2} breaks down when $p>2$.

\medskip

The preceding example involved a metric which, though $C^\infty$, is
not analytic.  It is also possible to show that \eqref{2.2} may break
down for a given $p>2$ when $n=3$ even if one considers analytic
metrics.

To see this we now let
\begin{equation}\label{alphak}
\alpha(s)=\alpha_k(s)=s^k, \, \, \, k=1,2,\dots .
\end{equation}
We then, for small $x$, let $\sum g_{jk}dx_jdx_k$ be the metric whose cometric is
$d\xi^2+2\alpha_k(x_2)d\xi_1d\xi_3$.  It then follows that for $x_1\in
\Bbb R$ and $-\pi<\theta<\pi$
\begin{equation}\label{pig}
t\to x(x_1,\theta;t)=(\, x_1+t\sin\theta,t\cos\theta,
\frac1{k+1}\sin\theta \cos^k\theta \, t^{k+1})\end{equation}
are geodesics.  Moreover, if we fix a small $x_2<0$, the last part of
the lemma ensures that we can find a small ball $B$ centered at
$(0,x_2,0)$ so that if $x\in B$ there is a unique geodesic as in
\eqref{pig} which passes through $x$.  Since $|t^{k+1}|<\delta$ if
$|t|<\delta^{1/(k+1)}$, if we fix $c>0$ and now let
$$f_\delta(x)=1\, \, \text{if } \, 0\le x_2\le \delta^{1/(k+1)}, \,
|x_1|\le c, \, |x_3|\le \delta, \, \, \text{and } \, f_\delta(x)=0\, \,
\text{otherwise},$$
then, if the center of $B$ is close to the origin,
$$\Mdel f_\delta (x)\ge c_0\delta^{1/(k+1)}, \, \, x\in B,$$
for some $c_0>0$ depending on $c$ and $B$.  Consequently,
$$\|\Mdel f_\delta\|_{L^1(B)}\, /\, \|f_\delta\|_{L^p}\ge
c_0'\delta^{1/(k+1)-(k+2)/(k+1)p}.$$ 
Since 
$$1-3/p>(k+2)/(k+1)p-1/(k+1) \, \, \text{when } \, p>(2k+1)/k,$$
it follows that \eqref{2.2} breaks down for a given fixed $p$ if $k$ is large.

\noindent{\bf Remark.}
Notice that when $k=1$ we only recover the trivial requirement for
\eqref{2.2} that $p\ge 3$.  To explain the difference between this case
and the others we note that in all cases, the key point involved the
behavior of the geodesics in the $(x_2,x_3)$ direction.  This is
dictated by the $R^3_{232}$ component of the curvature tensor.  A
calculation shows that, when $k=1$, $R^3_{232}=-(3-5x^2_2)/4(1-x^2_2)$,
and so in particular $R^3_{232}\approx -1/4$ when $|x_2|$ is small.  In
the other cases, where $k=2,3,\dots$, though, $R^3_{232}\approx
-x_2^{2k-2}$ near $x_2=0$ and so this sectional curvature vanishes to
higher and higher order at $x_2=0$ as $k\to +\infty$.  In the first
example of course it vanishes of infinite order.  Based on this and
related results to follow one might conjecture that for curved spaces
one would want to assume that the sectional curvatures are pinched away
from zero to obtain favorable bounds for Nikodym maximal
operators or related oscillatory integral operators.  This condition by
itself is probably not sufficient since even though the results of
\cite{W} seem to easily extend to the hyperbolic space setting, it
seems that the arguments in this paper can be used to show that
\eqref{cordconj} cannot hold for certain local perturbations of ${\Bbb
H}^n$ when $n$ is odd and $p>(n+1)/2$.

We hope to explore these points in a later work.

\newsection{Negative results for maximal operators in higher {\it odd}
dimensions}

It is not hard to adapt the argument for the three-dimensional case and
show that \eqref{2.2} does not hold in general for an odd-dimensional
Riemannian manifold when $(n+1)/2<p\le n$.  Later we shall see that the
inequality does hold though in the complimentary range where $1\le p\le
(n+1)/2$.  We shall then use this fact to show how, at least for odd
dimensions, our constructions give the maximum possible amount of
``focusing" of geodesics.

To prove the negative results for \eqref{2.2} when $n$ is odd we shall
consider cometrics on $T^*\Rn$ of the form
\begin{equation}\label{3.1}
\sum_{j,k=1}^n g^{jk}(x)d\xi_jd\xi_k=
d\xi^2+2\alpha(x_{(n+1)/2})\sum_{j=1}^{(n-1)/2}
d\xi_{(n+1)/2-j}d\xi_{(n+1)/2+j},
\end{equation}
where $\alpha\in C^\infty$ satisfies $|\alpha|<1$ and $\alpha(0)=0$.  We
then, as before, let $\sum g_{jk}(x)dx_jdx_k$ be the associated
Riemannian metric where $g_{jk}=(g^{jk})^{-1}$.  We then can use the
proof of Lemma \ref{mainlemma} to see that if
$\theta=(\theta_1,\dots,\theta_{(n-1)/2})$ is fixed and satisfies
$|\theta|^2=\sum \theta^2_j<1/2$, say, and if $(x_1,\dots,x_{(n-1)/2})$
is fixed, then
\begin{multline}
t\to
x(x_1,\dots,x_{(n-1)/2},\theta;t)
\\
=(x_1+t\theta_1,\dots,x_{(n-1)/2}+t\theta_{(n-1)/2},
t\sqrt{1-|\theta|^2},\theta\alpha^{(-1)}(t\sqrt{1-|\theta|^2})/\sqrt{1-|\theta|^2})
\end{multline}
parameterizes a geodesic.  As before $\alpha^{(-1)}$ denotes the
primitive of $\alpha$ vanishing at the origin.

In what follows we shall assume that $\alpha$ is given by
\eqref{alpha}.  Then our metric of course agrees with the Euclidean one
when $x_{(n+1)/2}\ge 0$.

Note that the Jacobian of the map sending
$$(x_1,\dots,x_{(n-1)/2},\theta,t)\to x(x_1,\dots,x_{(n-1)/2},\theta;t)$$
equals $|\alpha^{(-1)}(t)|^{(n-1)/2}$ when $\theta=0$.  Consequently,
if we fix $x_{(n+1)/2}<0$ we can find a ball $B$ centered at
$(0,\dots,0,x_{(n+1)/2},0,\dots,0)$ so that if $x\in B$ then there is a
unique geodesic $\gamma_x$ which contains $x$ and lies in the
$(n+1)/2$-plane $\Pi=\{x:\, x_j=0, \, (n+1)/2<j\le n\}$ when
$x_{(n+1)/2}>0$.  Consequently, if we assume, depending on our
definition of $\Mdel$, that the center of $B$ is sufficiently close to
the origin, we obtain
$$\Mdel f_\delta(x)\ge c_0>0, \, \, \, x\in B,$$
if for a given fixed $c>0$
$$
f_\delta(x)=
\begin{cases}
1\, \, \, \text{if } \, |(x_1,\dots,x_{(n+1)/2})|<c, \, \, \text{and }
\, 
|x_j|<\delta, \, \, (n+1)/2<j\le n
\\
0\, \, \, \text{otherwise}.
\end{cases}$$
From this we conclude that, for some $c'_0>0$, 
$$\|\Mdel f_\delta\|_{L^1(B)}\, /\, \|f_\delta\|_{L^p}\ge
c_0'\delta^{-(n-1)/2p}.$$
Since
$$n/p-1<(n-1)/2p\, \ \, \text{when } \, p>(n+1)/2,$$
we conclude that \eqref{2.2} cannot hold here for $p>(n+1)/2$.

This example of course involved a smooth metric which was not real
analytic.  As in the three-dimensional case, though, it is
straightforward to modify the construction using \eqref{alphak} to see
that given $p_0>(n+1)/2$ there is a real analytic metric for which
\eqref{2.2} cannot hold when $p_0<p\le n$.

\newsection{Negative results for maximal operators in higher even
dimensions}

The negative results for even dimensions are somewhat different since
we cannot have sharp focusing of space filling geodesics into an
$(n+1)/2$-dimensional submanifold since $(n+1)/2$ is not an integer
when $n$ is even.  In the next section we shall say a bit more about
the difference between even and odd dimensions.  In particular we shall
show that for $n$ even there can only be sharp focusing of space
filling geodesics into submanifolds of dimension $(n+2)/2$ when $n$ is
even.  Because of this fact our methods only show that \eqref{2.2}
cannot hold in general for $p>(n+2)/2$ on even dimensional curved
manifolds.

To prove this we shall consider cometrics of the form
\begin{equation}\label{4.1}
\sum_{j,k=1}^n
g^{jk}(x)d\xi_jd\xi_k=d\xi^2+2\alpha(x_{(n+2)/2})
\sum_{j=1}^{(n-2)/2}d\xi_{n/2-j}d\xi_{(n+2)/2+j},
\end{equation}
assuming as usual that $\alpha$ is smooth and that $|\alpha|<1$.  If
then $\sum g_{jk}(x)dx_jdx_k$ is the corresponding metric, one checks
using the earlier arguments that, when $(x_1,\dots,x_{n/2})$ and 
$\theta=(\theta_1,\dots,\theta_{(n-2)/2})$ with $|\theta|<1/2$ are
fixed, the curves
\begin{multline*}
t\to x(x_1,\dots,x_{n/2},\theta;t)
\\
=(x_1+t\theta_1,\dots,x_{(n-2)/2}+t\theta_{(n-2)/2},x_{n/2},t\sqrt{1-|\theta|^2},
\theta\alpha^{(-1)}(t\sqrt{1-|\theta|^2})/\sqrt{1-|\theta|^2})
\end{multline*}
are geodesic.

If we assume that $\alpha$ is as in \eqref{alpha} then the Jacobian of
$$(x_1,\dots,x_{n/2},\theta,t)
\to 
x(x_1,\dots,x_{n/2},\theta;t)$$
is nonsingular when $\theta=0$ and $t<0$.  Consequently, if we fix
$x_{(n+2)/2}<0$ and $x_{n/2}\in \Bbb R$ there is a ball $B$ centered at
$(0,\dots,x_{n/2},x_{(n+2)/2},0,\dots,0)$ so that if $x\in B$ there is
a unique geodesic $\gamma_x$ containing $x$ and lying in the
$(n+2)/2$-plane $\Pi=\{x: \, x_j=0, \, (n+2)/2<j\le n\}$ when
$x_{(n+2)/2}\ge 0$.

To use this, for a given $c>0$, we put
$$f_\delta(x)=
\cases
1\, \, \, \text{if } \, \, |(x_1,\dots,x_{(n+2)/2})|<c, \, \, \text{and
} \, |x_j|<\delta, \, \, (n+2)/2<j\le n
\\
0\, \, \, \text{otherwise}.
\endcases$$
Then if the center of $B$ is close to the origin, we must as before
have that $\Mdel f_\delta(x)$ is bounded below by a positive constant
(depending on $B$) for each $x\in B$.  We then conclude that, for some
$c_0>0$,
$$\|\Mdel f_\delta\|_{L^1(B)}\, / \, \|f\|_{L^p} \ge
c_0\delta^{-(n-2)/2p},$$
which implies that \eqref{2.2} cannot hold for $p>(n+2)/2$ since
$n/p-1<(n-2)/2p$ for such $p$.

\newsection{Bounds for maximal functions and lower bounds on the
dimension of 
\\
Nikodym-type sets}

The main result of this section is the following

\bth\label{posthm}  Let $(M^n,g)$ be a complete Riemannian manifold of
dimension $n\ge2$, and let $\Mdel$ be as in \eqref{2.1} where $r=\min
\{1, (\text{inj }M^n)/2\}$, with $\text{inj }M^n$ denoting the
injectivity radius of $M^n$.  If then $K\subset M^n$ is a fixed compact
set
\begin{multline}\label{pos}
\|\Mdel f\|_{L^q(K)}\le
C_{p,\varepsilon}\delta^{1-n/p-\epsilon}\|f\|_{L^p},
\\
\text{if } \, \, \supp f\subset K, \, \, 1\le p\le (n+1)/2 \, \,
\text{and } \, q=(n-1)p'.
\end{multline}
\eth

In view of our earlier negative results \eqref{pos} is best possible in
the general curved space setting when $n$ is odd.

Before turning to the proof, let us see how \eqref{pos} and our earlier
constructions yield sharp lower bounds for the dimension of
Nikodym-type subsets of general odd-dimensional 
manifolds.\footnote{The sets actually correspond to sets which in the Euclidean setting
would contain compliments of the usual Nikodym sets (see \cite{Fal});
however, we are following the terminology in \cite{B1}.
}

\noindent{\bf Definition.}  If $\Pi\subset \subset M^n$ let $\Pi^*$ denote all
points $x\in M^n$ for which there is a geodesic $\gamma_x\ni x$ of
length $\le r=\min \{1, (\text{inj }M^n)/2\}$ which intersects $\Pi$ in
a set of positive length, that is, $|\Pi\cap \gamma_x|>0$.  We then
call $\Pi$ a {\it Nikodym-type set} if $\Pi^*$ has positive measure.

\bco\label{nikbound}  If $\Pi$ is a Nikodym-type subset of $M^n$ then
the Minkowski dimension of $\Pi$ is at least $(n+1)/2$.
\eco

For odd $n$ the lower bounds are sharp since we have shown that if the
cometric is as in \eqref{3.1} with $\alpha$ given by \eqref{alpha},
then the intersection of the $(n+1)/2$-plane $\{x:\, x_j=0, \,
(n+1)/2<j\le n\}$ with any ball centered at the origin is a
Nikodym-type set.  Also, the corollary implies that if $\Pi$ is a
submanifold and a Nikodym-type set then its dimension must be $(n+2)/2$
for even $n$.  This accounts for the difference between our negative
results in even and odd dimensions since our strongest counterexamples
all involve such sets.

The proof of the corollary is very simple.  We must show that if $\Pi$
is a Nikodym-type set then for every $\epsilon>0$ there is a constant
$c_\epsilon>0$ so that
\begin{equation}\label{5.2}
|\Pi^\delta|\ge c_\epsilon\delta^{(n-1)/2+\epsilon}, \, \, 0<\delta\le
1
\end{equation}
if $\Pi^\delta$ denotes a $\delta$-neighborhood of $\Pi^\delta$.
To show this we simply note that
$$\Pi^* \subset \cup_{\lambda>0}\, \{x:\, \inf_{0<\delta\le1}(\Mdel
\chi_{\Pi^\delta})(x)>\lambda\, \}$$
if $\chi_{\Pi^\delta}$ denotes the characteristic function of
$\Pi^\delta$.  Hence, if $\lambda>0$ is small and fixed
$$|\{x:\, \inf_{0<\delta\le 1}(\Mdel \chi_{\Pi^\delta})(x)>\lambda \,
\}| \ge c_0>0$$
if $|\Pi^*|>0$.  Since $\lambda$ is fixed, we conclude from \eqref{pos}
with $p=(n+1)/2$ (see also \eqref{weak} below) that if $\epsilon>0$ 
$$0<c'_0\le C_{\lambda,\epsilon}\delta^{1-n-\epsilon}|\Pi^\delta|^2, \,
\, 0<\delta\le 1,$$
which of course yields \eqref{5.2} and completes the proof.

\medskip

Turning to the proof of Theorem \ref{posthm}, let us first point out that 
undoubtedly one does not have to assume, in the definition of
$\Mdel$, that $|\gamma_x|$ is smaller
than a multiple of the injectivity radius (cf. \cite{S3}), but one
needs this hypothesis to be able to use the simple arguments of
Bourgain \cite{B1} and Wolff \cite{W}.  To see where this restriction is used we need to
introduce some notation.  If $\gamma_j(s)$, $s\in [\alpha_j,\beta_j]$
are two geodesics parameterized by arclength we set
$$\theta(\gamma_1,\gamma_2)=\min_{s_j\in[\alpha_j,\beta_j]}
\text{dist}((x_1(s_1), x'_1(s_1)),(x_2(s_2),x'(s_2))).$$
Here $\text{dist }$ comes from the natural metric on the unit cosphere
bundle induced by our given Riemannian metric on $M^n$.  Also, if $a\in
M^n$ and $\lambda>0$ let $B(a,\lambda)$ denote the geodesic ball 
radius $\lambda$ centered at $a$.

With this notation we shall require the following simple result which
is essentially contained in \cite{S1}.

\ble\label{spread}  Suppose that $\gamma_j$, $j=1,2$ are geodesics whose
length does not exceed $r=\min\{1, (\text{inj }M^n)/2\}$ and which belong 
to a fixed compact subset $K\subset
M^n$.  Suppose also that $a\in
T^\delta_{\gamma_1}\cap T^\delta_{\gamma_2}$.  Then there is a constant
$c>0$, depending on $(M^n,g)$ and $K$, but not on $\delta>0$ and
$0<\lambda\le 1$, so that
\begin{equation*}
(T^\delta_{\gamma_1}\cap T^\delta_{\gamma_2})\backslash
B(a,\lambda)=\emptyset \quad \text{if } \, \,
\theta(\gamma_1,\gamma_2)\ge \delta/c\lambda.
\end{equation*}
\ele

To proceed, we need to make a couple of easy reductions.  
We first
notice that since we are assuming that $\supp f\subset K$, where $K$ is
a fixed compact subset of $M^n$, it suffices to show that the variant of
\eqref{pos} holds where in the left side the norm is taken over a fixed
compact subset of a coordinate patch.  
We can even assume further, for
the sake of convenience, that local coordinates have been chosen so
that the vertical lines  where $x'=(x_1,\dots,x_{n-1})$ is constant are
all geodesic.  It then suffices to show that, if in our definition of
$\Mdel$ we add the restriction that $\gamma_x$ satisfies
$\theta(\gamma_x,\ell)\le c_0$ for some such line $\ell$ and
a given small constant $c_0>0$, then \eqref{pos}
holds.   This in turn would be a
consequence of the stronger bounds
$$\bigl(\, \int |\Mdel f(x')|^q\, dx'\, \bigr)^{1/q}\le
C_\epsilon
\delta^{1-n/p-\epsilon}\|f\|_p, \quad q=(n-1)p/(p-1), \, \, 1\le p\le
(n+1)/2,$$
assuming as before that $f$ has small support, and that now
$$\Mdel f(x')=\Mdel f(x',0).$$
Here and in what follows we are assuming that $x'\in K'=\{x\in K: \, x_n=0\}$.

Since the bound for $p=1$ is trivial, the preceding inequality would follow from
showing that, under the above assumptions, the maximal operator is of
restricted weak-type $((n+1)/2,n+1)$ with norm $O(\delta^{(1-n)/(n+1)})$.  To be
more specific, we need to show that if $E$ is contained in a fixed
compact subset of a coordinate patch as above then
\begin{equation}\label{weak}
|\{x': \, \Mdel \chi_E(x')>\lambda\}|\le
C\lambda^{-(n+1)}\delta^{1-n}|E|^2.
\end{equation}
Since the set in question is empty for $\lambda>1$ we need only
consider $0<\lambda\le 1$.  To simplify the notation and arguments to
follow, we shall also let  $A$  denote a fixed large constant which is to be specified
later that depends
on $(M^n,g)$ and our support assumptions.
It then suffices to verify that
\begin{equation}\label{9'}
|\{x': \, \Mdel \chi_E(x')>A\lambda \}|\le
C\lambda^{-(n+1)}\delta^{1-n}|E|^2, \, \, \delta, \lambda\in (0,1],
\end{equation}
with $C$ here being equal to $A^{-(n+1)}$ times the constant in the
preceding inequality.

Assuming that $A$ is as above we choose a maximally
$A\delta/\lambda$-separated subset
$$\{x'_j\}_{j=1}^M={\cal I}$$
in $\{x': \, \Mdel \chi_E(x')>A\lambda\}.$ If we then note that
\begin{equation}\label{10}
|\{x': \Mdel \chi_E(x')>A\lambda\}|\le CM \cdot
(A\delta/\lambda)^{n-1},
\end{equation}
we conclude that our task is equivalent to obtaining an appropriate
upperbound on the cardinality $M$ of ${\cal I}$.

The first step in doing this is to notice that given $x_j'\in {\cal I}$
we can choose a geodesic $\gamma_j$ containing $(x',0)$ of length $\le
r$ so that
\begin{equation}\label{11}
|E\cap T^\delta_{\gamma_j}|\ge A\lambda |T^\delta_{\gamma_j}|.
\end{equation}
Since $|T^\delta_{\gamma_j}|\approx \delta^{n-1}$, if we sum over $j$,
we conclude that
$$\sum_{j=1}^M|E\cap T^\delta_{\gamma_j}|\ge c_0M\lambda \delta^{n-1}$$
for a fixed constant $c_0>0$.

From this we conclude that there must be a point $a\in E$ belonging to
at least
$$N=c_0M\lambda\delta^{n-1}/|E|$$ 
of the tubes $T^\delta_{\gamma_j}$.  Label these as
$\{T^\delta_{\gamma_{j_k}}\}_{1\le k\le N}$.

If we invoke the preceding lemma, we conclude that
$(T^\delta_{\gamma_{j_1}}\cap T^\delta_{\gamma_{j_2}})\backslash
B(a,\lambda)=\emptyset$ if $\theta(\gamma_{j_1},\gamma_{j_2})\ge
\delta/c\lambda$, with $c>0$ being a fixed constant.  Since ${\cal I}$ is
$A\delta/\lambda$-separated, 
this condition is automatically
satisfied for $j_1\ne j_2$ if $A$ is large enough, assuming, as above, that the
geodesics are close to vertical lines.
This in turn
implies that the tips of the tubes
$\tau^\delta_{j_k}=T^\delta_{\gamma_{j_k}}\backslash B(a,\lambda)$,
$1\le k\le N$, are disjoint.  Since
$$|T^\delta_{\gamma_j}\cap B(a,\lambda)|\le
C_0\lambda|T^\delta_{\gamma_j}|$$
for a fixed constant $C_0$, we conclude from \eqref{11} that if we also
assume that $A\ge 2C_0$, then
$$|\tau^\delta_{\gamma_{j_k}}\cap E|\ge A\lambda
|T^\delta_{\gamma_{j_k}}|/2,
\quad 1\le j\le N.$$
Hence, if we sum and use the aforementioned disjointness, we conclude
that
$$|E|\ge \sum_{j=1}^N|\tau^\delta_{j_k}\cap E|\ge AN\lambda
\delta^{n-1}/2\ge CM\lambda^2\delta^{2(n-1)}/|E|.$$
Since this yields
$$M\le C'\lambda^{-2}\delta^{-2(n-1)}|E|^2,$$
we obtain \eqref{9'} from \eqref{10}, which completes our proof.

\newsection{Negative results for oscillatory integrals in odd
dimensions}

In the remainder of the paper we shall show that bounds of the form
\eqref{oscest} need not hold for certain $2n/(n-1)<q<2(n+1)/(n-1)$ if
$n>2$ and
\begin{equation}\label{6.1}
(S_\lambda f)(x)=\int e^{i\lambda \dist(x,y)} a(x,y) f(y)\, dy,
\end{equation}
with $\dist(x,y)$ denoting the Riemannian distance between $x$ and $y$
in $\Rn$ measured by a non-Euclidean metric.  To avoid the singularity
of the phase we shall assume that $a$ vanishes near the diagonal and
for convenience we shall also assume that $0\le a\in
C^\infty_0(\Rn\times \Rn)$ and that
\begin{equation}\label{6.2}
a(x,y)\ne 0\, \, \, \text{if } \, x=0\, \, \text{and } \, \, y_j=0, \,
j\ne (n+1)/2, \, \, y_{(n+1)/2}=-1.
\end{equation}

Here we are assuming that $n\ge3$ is odd.  We then take our metric to
be dual to the one in \eqref{3.1} where $\alpha$ is given
\eqref{alpha}.

To proceed, we need to use an argument from Bourgain \cite{B1}.  (See
also Fefferman \cite{F1}.)  To be more specific, we first need to
recall that if, for every $\varepsilon>0$, 
$S_\lambda: L^p\to L^q$ with norm $C_{p,q}\le
C_\varepsilon\lambda^{-n/q+\varepsilon}$, then the adjoint operator
\begin{equation}\label{6.3}
(S^*_\lambda g)(y)=\int e^{-i\lambda \dist(x,y)}a(x,y)g(x)\, dx
\end{equation}
must send $L^{q'}\to L^{p'}$ with the same norm.  Finally, we need to
recall (see p. 484, Theorem 2.7 in \cite{Rubbook} or \cite{St1}) 
that the dual bounds in turn imply a vector
valued version
\begin{equation}\label{6.4}
\|\, (\sum_\alpha|S^*_\lambda g_\alpha|^2)^{1/2}\, \|_{p'}\le
C_\varepsilon'\lambda^{-n/q+\varepsilon} 
\|\, (\sum_\alpha|g_\alpha|^2)^{1/2}\, \|_{q'}, \quad \varepsilon>0,
\end{equation}
with $C_\varepsilon'$ being a fixed multiple of $C_\varepsilon$ 
for a given $p$ and $q$.

To show that this inequality need not hold for certain $q>2n/(n-1)$,
let $y$ be as in \eqref{6.2}.  We then can find a ball $B$ centered at
$y$ so that if $z\in B$ there is a unique geodesic $\gamma_z\ni z$ which
is contained in the $(n+1)/2$-plane $\{x:\, x_j=0, \, (n+1)/2<j\le n\,
\}$ when $x_{(n+1)/2}\ge 0$.  We then choose a maximally
$\lambda^{-1/2}$-separated set of points $z_\alpha\in B\cap \{y: \,
y_{(n+1)/2}=-1\}$.   We also define the Euclidean cylinders
\begin{equation}\label{6.5}
T_\alpha =\{x:\, x_{(n+1)/2}\ge 0, \, \, |x|\le 1, \, \,
\dist(x,\gamma_{z_\alpha})\le c\lambda^{-1/2}\},
\end{equation}
and set
$$g_\alpha(x)=e^{i\lambda\dist(x,z_\alpha)}\chi_{T_\alpha}(x).$$
Keeping \eqref{6.2} in mind, if $c>0$ in \eqref{6.5} and the diameter of
$B$ are small enough, one checks that
$$|S^*_\lambda g_\alpha(y)|\approx |T_\alpha|\approx
\lambda^{-(n-1)/2}, \, \, \text{if } \, \dist(y,\gamma_{z_\alpha})<c\lambda^{-1/2} \,
\, \text{and } \, y\in B,$$
using the fact that $\nabla_x(\, \dist(x,z_\alpha)-\dist(x,y)\, )=0$ if
$x,y\in \gamma_{z_\alpha}$.
Thus,
\begin{equation}\label{6.6}
\lambda^{-(n-1)/2}\approx \int_{B} \max_\alpha |S^*_\lambda
g_\alpha(y)|\, dy \le \int_B (\sum_\alpha |S^*_\lambda
g_\alpha|^2)^{1/2}\, dy.
\end{equation}
If we use H\"older's inequality and \eqref{6.4} we can dominate the
right hand side by
\begin{equation}\label{6.7}
C_\varepsilon\lambda^{-n/q+\varepsilon}
\|\, (\sum |g_\alpha|^2)^{1/2}\, \|_{q'}=C_\varepsilon
\lambda^{-n/q+\varepsilon}
\|\, (\sum\chi_{T_\alpha})^{1/2}\, \|_{q'}.
\end{equation}

Recall that $\chi_{T_\alpha}(x)=0$ outside of the intersection of the
unit ball with the slab where $|x_j|\le c\lambda^{-1/2}$, $(n+1)/2<j\le
n$ and $x_{(n+1)/2}\ge0$.  In this region the metric is Euclidean and
it is not hard to see by a simple volume packing argument that a given
point $x$ in the region can lie in at most $O(\lambda^{(n-1)/4})$ of
the cylinders $T_\alpha$.  This just follows from the fact that there
are $O(\lambda^{(n-1)/2})$ cylinders of volume $\approx
\lambda^{-(n-1)/2}$ uniformly distributed in the above set which has volume
$\approx \lambda^{-(n-1)/4}$.  

If we use this overlapping bound, we conclude that
\begin{equation}\label{6.8}
\|(\sum_\alpha \chi_{T_\alpha})^{1/2}\|_{q'}\le
C\lambda^{(n-1)/8}\lambda^{-(n-1)/4q'}.
\end{equation}
If we combine this with the preceding two inequalities we conclude that
if the equivalent version \eqref{6.4} of \eqref{oscest} held, then 
as $\lambda\to +\infty$ we would have
$$\lambda^{-(n-1)/2}\le
C_\varepsilon\lambda^{-n/q+\varepsilon}
\lambda^{(n-1)/8}\lambda^{-(n-1)/4q'}, \quad \forall \varepsilon>0.$$ 
This in turn leads to the condition that
$$q\ge q_n=2(3n+1)/3(n-1)>2n/(n-1)$$
even if the weaker version,
$$\|S_\lambda f\|_{q}\le
C_\varepsilon\lambda^{-n/q+\varepsilon}\|f\|_\infty,
\quad \varepsilon>0,$$
of \eqref{oscest} held.  In particular, we conclude that when $n=3$
\eqref{oscest} breaks down in the curved space setting for $3\le
q<10/3$.
Also, as before, one could modify this construction and show that for a
given $2n/(n-1)<q<q_n$ \eqref{oscest} need not hold even on a manifold
with an analytic metric.

\newsection{Negative results for oscillatory integrals in even higher
dimensions}

It is easy to adapt the above argument and show that \eqref{oscest}
need not hold for certain $2n/(n-1)<q<2(n+1)/(n-1)$ when $n\ge 4$ is
even.  One lets the Riemannian metric on $\Rn$ correspond to the
cometric \eqref{4.1} where, as before, $\alpha$ is as in \eqref{alpha}.

One then replaces \eqref{6.2} with the condition that $a(x,y)\ne 0$
when $x=0$ and $y_j=0$, $j\ne (n+2)/2$, and $y_{(n+2)/2}=-1$.  One
makes similar modifications of the other parts of the proof for odd
$n$, replacing $(n+1)/2$ by $(n+2)/2$.  Then \eqref{6.6} and
\eqref{6.7} go through.  Inequality \eqref{6.8}, though, must be
modified since the cylinders $T_\alpha$ now lie in the slab where
$|x_j|\le c\lambda^{-1/2}$, $(n+2)/2<j\le n$, $x_{(n+2)/2}\ge 0$ and
$|x|\le 1$.  The arguments for the odd-dimensional case imply that a
point in this region belongs to $O(\lambda^{(n-2)/4})$ of the
$T_\alpha$.  Consequently, \eqref{6.8} must be replaced in even
dimensions by
$$\|(\sum_\alpha \chi_{T_\alpha})^{1/2}\|_{q'}\le
C\lambda^{(n-2)/8}\lambda^{-(n-2)/4q'}.$$
If we combine this with \eqref{6.6} and \eqref{6.7} we conclude that
if \eqref{oscest} holds for this example then we must have
$$\lambda^{-(n-1)/2}\le
C_\varepsilon\lambda^{-n/q+\varepsilon}\lambda^{(n-2)/8}\lambda^{-(n-2)/4q'},
\quad \forall \varepsilon>0,$$
as $\lambda\to +\infty$.  This in turn leads to the condition that for
even $n\ge4$ we must have $q\ge 2(3n+2)/(3n-2)>2n/(n-1)$.

\end{document}